\documentclass[jsco]{academic}

\usepackage{xy}

\input xy
\xyoption{all}

\usepackage{amsmath}
\usepackage{amssymb}
\usepackage{amsfonts}

\def\bC{{\mathbb C}}

\def\bP{{\mathbb P}}

\def\cO{{\mathcal O}}

\def\lar{\longrightarrow}

\def\iso{{\, \cong\, }}

\def\<{\langle}
\def\>{\rangle}

\newtheorem{example}[theorem]{Example}

\begin{document}
\shortauthor{Shaska}
\shorttitle{Decomposable Jacobians} 
\title{Curves of genus 2 with (N,N) \\
decomposable Jacobians}
\author{T. Shaska}

\address{Department of Mathematics, University of Florida, Gainesville, FL 32611.}


\maketitle

\begin{abstract} Let $C$ be a curve of genus 2 and $\psi_1:C \lar E_1$ a map
of degree $n$, from $C$ to an elliptic curve $E_1$, both curves defined over $\bC$.  
This map induces a degree
$n$ map $\phi_1:\bP^1 \lar \bP^1$ which we call a Frey-Kani covering.  
We determine all possible
ramifications for $\phi_1$. If $\psi_1:C \lar E_1$ is maximal then
there exists a maximal map $\psi_2:C\lar E_2$, of degree $n$, to some
elliptic curve $E_2$ such that there is an isogeny of degree $n^2$ from 
the Jacobian $J_C$ to $E_1 \times E_2$.
We say that  $J_C$ is $(n,n)$-decomposable. If the degree $n$ is odd the pair $(\psi_2,
E_2)$ is canonically determined.
For $n=3,\, 5$,  and 7,  we give arithmetic examples of curves whose
Jacobians are $(n,n)$-decomposable.
\end{abstract} 
%*********************************************************************************
\section{Introduction}
\par Curves of genus 2 with non-simple Jacobians are of much
interest. Their 
Jacobians have large torsion subgroups,  e.g. 
Howe, Lepr\'evost, and Poonen
 have found a family of genus 2 curve with 128 rational points in its
Jacobian, see \cite{P}.  For other   
applications of genus 2 curves with $(n,n)$-decomposable Jacobians see
Frey \cite{Fr}. In this paper, we 
discuss  genus 2 curves $C$  whose function fields have  maximal elliptic
subfields. These  elliptic subfields occur in pairs $(E_1, E_2)$ and
we call each 
the complement of the other in $J_C$. The Jacobian of
$C$ is isogenous to $E_1\times E_2$.  Let  $\psi:
C\to E$ be a maximal cover (cf. section 4) of odd degree $n$.  The
moduli space parameterizing these covers is a surface, more precisely 
the  product of modular
curves $X(n)\times X(n)/\Delta$, see  Kani \cite{K}.  When $\psi:
C\to E$ is degenerate (cf. section 2), this moduli space is a
curve. Getting algebraic descriptions for these spaces is extremely
difficult for large $n$ (e.g. $n\geq 7$).  Also, one would like to know how the elements
of the pair $(E_1, E_2)$ relate to each other. 

\par In sections 2 and 3 we define a Frey-Kani covering and determine
all their possible ramifications. In section 4 we consider maximal
covers. These covers allow us to determine the complement of $E_1$
uniquely. The last section deals with some applications when
$n=3,5$, or 7. 
%*************************************************************************
\section{Frey - Kani covers}
%*************************************************************************
\par Let $C$ and $E$ be curves of genus 2 and 1, respectively. Both are smooth, 
projective curves defined over $\bC$.  Let 
$\psi: C \longrightarrow E$ be a covering of degree $n$.   We say that $E$ is an degree $n$ 
{\it elliptic subcover } of $C$.  From the
Riemann-Hurwitz formula, $\sum_{P \in C}\, (e_{\psi}\,(P) -1)=2$ where
$e_{\psi}(P)$ is the ramification index of points $P \in C$, under
$\psi$.  Thus,  we have two points of ramification index 2 or one point of
ramification index 3.  The two points of ramification index 2 can be
in the same fiber or in different fibers.  Therefore,  we have the following
cases of the covering $\psi$:

{\bf Case I.} There are $P_1$, $P_2 \in C$, such that
$e_{\psi}({P_1})=e_{\psi}({P_2})=2$, $\psi(P_1) \neq \psi(P_2)$, and
	$\forall  P \in C\setminus \{P_1,P_2\}$,  $e_{\psi}(P)=1$.

{\bf Case II.} There are $P_1$, $P_2 \in C$, such that
$e_{\psi}({P_1})=e_{\psi}({P_2})=2$, $\psi(P_1) = \psi(P_2)$, and
	$\forall  P \in C\setminus \{P_1,P_2\}$,  $e_{\psi}(P)=1$.

{\bf Case III.} There is $P_1 \in C$ such that $e_{\psi}(P_1)=3$, and
$ \forall P \in C \setminus \{P_1\}$, $e_{\psi}(P)=1$

In case I (resp. II, III) the cover $\psi$ has 2 (resp. 1) branch
points in E.

\par Denote the hyperelliptic involution of $C$ by $w$. We choose
$\mathcal O$ in E such that $w$ restricted to $E$ is the hyperelliptic
involution on E, see \cite{FK} or \cite{Ku}. We denote the restriction
of $w$ on $E$ by $v$, $v(P)=-P$.  Thus, $\psi \circ w=v \circ
\psi$. E[2] denotes the group of 2-torsion points of the elliptic
curve E, which are the points fixed by $v$.   
The proof of the following two lemmas is straightforward and will be omitted.

\begin{lemma} \label{lem1}
a) If $Q \in E$, then $\forall P \in \psi^{-1}(Q)$,
$w(P) \in \psi^{-1}(-Q)$.

b) For all $P\in C$, $e_\psi(P)=e_\psi\,({w(P)})$.  
\end{lemma} 

\par Let $W$ be the set of points in C fixed by $w$.  Every curve of
genus 2 is given, up to isomorphism, by a binary sextic, so there are
6 points fixed by the hyperelliptic involution $w$, namely the Weierstrass points of 
$C$.  The following lemma 
determines the distribution of the Weierstrass points in fibers of 
2-torsion points.

\begin{lemma}\label{lem2}
\begin{enumerate}
\item $\psi(W)\subset E[2]$
\item If $n$ is an odd number then
\subitem i) $\psi(W)=E[2]$
\subitem ii) If $ Q \in E[2]$ then  \#$(\psi^{-1}(Q) \cap W)=1 \mod (2)$
\item If $n$ is an even number then for all $Q\in E[2]$,  
\#$(\psi^{-1}(Q) \cap W)=0 \mod (2)$
\end{enumerate}
\end{lemma}

\par Let $\pi_C: C \lar \bP^1$ and $\pi_E:E \lar \bP^1$ be the
natural degree 2 projections.  The hyperelliptic involution permutes
the points in the fibers of $\pi_C$ and $\pi_E$.  The  ramified
points of $\pi_C$, $\pi_E$ are respectively points in $W$ and $E[2]$
and their ramification index is 2.
There is  $\phi:\bP^1 \lar \bP^1$ such that 
the diagram commutes, see  Frey \cite{FK} or Kuhn \cite{Ku}.

$$\begin{matrix}
C & \buildrel{\pi_C}\over\lar & \bP^1\\
\psi \downarrow &  & \downarrow \phi \\
E & \buildrel{\pi_E}\over\lar & \bP^1
\end{matrix}$$

%$$\CD
%C @>{\pi_C}>> \bP^1 \\
%@V{\psi}VV        @VV{\phi}V \\
%E @>\pi_E>> \bP^1
%\endCD
%$$ 

\par The covering $\phi:\bP^1 \lar \bP^1$ will be called the corresponding
{\bf Frey-Kani covering} of $\psi:C \lar E$. It has first appeared in
\cite{FK} and \cite{Fr}. The term, Frey-Kani covering, has first been used by Fried in
\cite{Fried}.

\section{The ramification of Frey-Kani coverings}
%***********************************************************

\par In this section we  will determine the ramification of Frey-Kani coverings
$\phi:\bP^1 \lar \bP^1$. 
First we fix some notation. For a given branch point we will denote the ramification of
points in its fiber as follows. Any point $P$ of ramification index $m$  is denoted
by $(m)$.  If there are $k$ such points then we write
$(m)^k$.  We omit writing symbols for unramified points, in other
words $(1)^k$ will not be written.  Ramification data between two
branch points will be separated by commas. We denote by
 $\pi_E (E[2])=\{q_1, \dots , q_4\}$ and $\pi_C(W)=\{w_1,
\dots ,w_6\}$.

\subsection{The case when $n$ is odd}
%***************************************************************
The following theorem classifies the ramification types for the
Frey-Kani coverings $\phi:\bP^1 \lar \bP^1$ when the degree $n$ is odd.

\begin{theorem} \label{thm1}
If  $\psi:C \lar E$ is a covering of odd degree  $n$  
then the three cases of ramification
for $\psi$ induce the following cases for $\phi: \bP^1 \lar \bP^1$.
\begin{description}
\item[Case I:] (the generic case)

\subitem 	$\left ( (2)^\frac {n-1} 2  , (2)^\frac {n-1} 2  ,
(2)^\frac {n-1} 2  , (2)^ \frac {n-3} 2 ,(2)^1  \right )
$
\par Or  the following degenerate cases:

\item[Case II:] (the 4-cycle case and the dihedral case)

\subitem 	i) $ \left ((2)^\frac {n-1 }{ 2}  ,(2)^\frac{n-1 }{ 2}  , 
(2)^\frac {n-1 }{ 2} , (4)^1  (2)^\frac {n-7 }{ 2}    \right )$

\subitem	ii)	$\left ((2)^\frac {n-1}{ 2}  ,  (2)^\frac {n-1}{2}  ,  
      (2)^\frac {n-1}{ 2}, (2)^\frac {n-1} 2 \right )$

\subitem	iii)   $\left ((2)^\frac {n-1}{ 2}  ,  (2)^\frac {n-1}{2}  ,  
      (4)^1  (2)^\frac {n-5}{ 2}, (2)^\frac {n-3} 2 \right )$

\item[Case III:] (the 3-cycle case)

\subitem 	i) $ \left ((2)^\frac {n-1 }{ 2}  , (2)^\frac {n-1 }{ 2}  ,   
(2)^\frac {n-1}{ 2}  , (3)^1 (2)^\frac {n-5 }{ 2}  \right )$

\subitem 	ii) $\left ( (2)^\frac {n-1} 2  ,  (2)^\frac {n-1} 2  ,
(3)^1 (2)^\frac {n-3}  2,  (2)^\frac {n-3} 2
\right )$
\end{description}
\end{theorem}

\proof  From lemma \ref{lem2} we can assume that $\phi(w_i)=q_i$ for
$i\in \{1,2,3\}$ and $\phi(w_4)=\phi(w_5)=\phi(w_6)=q_4$.  Next we
consider the three cases for the ramification of $\psi:C\lar E$ and see
what ramifications they induce on $\phi:\bP^1\lar \bP^1$.

\par Suppose that $P\in \psi^{-1}(E[2])\setminus W $ and
$e_\psi(P)=1$. Then $e_{\psi}(P)\cdot e_{\pi_E}(\psi(P))=e_{\pi_C}(P) \cdot
e_{\phi}(\pi_C(P))=2$, so $e_{\phi}(\pi_C(P))=2$.

\par {\bf Case I:} There are $P_1$ and $P_2$ in $C$ such that
$e_{\psi}\, (P_1)=e_{\psi}\, (P_2)=2$ and $\psi(P_1)\neq\psi(P_2)$.  By
lemma \ref{lem1}, $e_\psi\,  ({w(P_1)})=2$.  So $w(P_1)=P_1$ or $w(P_1)=P_2$.

\par Suppose that $w(P_1)=P_1$, so $P_1 \in W$. If
$\pi_C(P_1)=w_i$ for $i \in \{1,2,3\}$, say $\pi_C(P_1)=w_1$, then
$e_{\pi_E \circ \psi}(P_1)=e_{\phi \circ \pi_C}(P_1)=4$, which implies
that $e_{\phi} (w_1)=2$.  All other points in the fiber of
$\pi_E\circ \psi (P_1)=:q_1$ have ramification index 2 under
$\phi$. So $\phi $ has even degree, which is a contradiction. If 
 $\pi_C (P_1)=w_i$ for $i \in \{4,5,6\}$, say $\pi_C (P_1)=w_4$, then in
the fiber of $q_4$ are: $w_4$ of ramification index 2, $w_5$ and $w_6$
unramified, and all other points have ramification index 2. So $\#
(\phi^{-1} (q_4) )=2+1+1+2k$, is even. 
Thus $P_1, P_2 \notin W$.  Then $P_1, P_2 \notin \psi^{-1}(E[2])$,
otherwise they would be in the same fiber.

\par Thus $P_2=w(P_1)\in C\setminus \psi^{-1}(E[2])$ and
$\psi(P_1)=-\psi(P_2)$.  Let $\pi_E \circ \psi (P_1)=\pi_E \circ \psi
(P_2)=q_5$ and $\pi_C(P_1)=\pi_C(P_2)=S$.  So $e_{\psi}(P_1)\cdot
e_{\pi_E}\, (\psi(P_1))=e_{\pi_C}\, (P_1) \cdot e_{\phi}\,(\pi_C(P_1))$.  Thus,
$e_{\phi}\, (\pi_C(P_1))=e_{\phi}\, (S)=2$.  All other points in
$\phi^{-1}(q_5)$ are unramified.

For $P \in W$, $e_{\pi_C}(P)=2$.  Thus $e_{\phi}\, (\pi_C(P))=1$.  All
$w_1, \dots w_6$ are unramified and other points in $\phi^{-1}(E[2])$
are of ramification index 2. By the Riemann - Hurwitz formula, $\phi$
is unramified everywhere else.  \par Thus, there are $\frac {n-1} 2$
points of ramification index 2 in the fibers $\phi^{-1}(q_1)$,
$\phi^{-1}(q_2)$, $\phi^{-1}(q_3)$, $\frac {n-3} 2$ points of
ramification index 2 in $\phi^{-1}(q_4)$, and one point of index 2 in
$\phi^{-1}(q_5)$.

\par {\bf Case II:} In this case, there are distinct $P_1$ and $P_2$
in $C$ such that $e_{\psi}(P_1)=e_{\psi}(P_2)=2$ and
$\psi(P_1)=\psi(P_2)$. Then $P_2=w(P_1)$ or $w(P_i)=P_i$, for $i=1,2$.

\par Let $P_1$ and $P_2$ be in the fiber which has three Weierstrass points. 

i) Suppose that $w$ permutes $P_1$ and $P_2$. So $P_1 $ and $P_2$ are
not Weierstrass points.  Then $e_{\pi_E \circ \psi}(P_1)=e_{\psi}(P_1)
\cdot e_{\pi_E}(\psi(P_1))=4$.  Thus $e_{\pi_C}(P_1)\cdot e_{\phi}(\pi_C
(P_1))=4$.  Since $e_{\pi_C}(P_1)=1$ then
 $e_{\phi}(\pi_C (P_1)=4$. So there is a point of index 4 in the fiber
of $q_4$. The rest of the points are of ramification index 2, as in
previous case, other then the $w_1, \dots , w_6$ which are unramified.

ii) Suppose that $w$ fixes $P_1$ and $P_2$. Thus $P_1$ and $P_2$ are
Weierstrass points.  Then $e_{\psi}(P_i)\cdot
e_{\pi_E}(\psi(P_i))=e_{\pi_C}(P_i) \cdot e_{\phi}(\pi_C(P_i))=4$.  So
$e_{\phi}(\pi_C(P_i))=2$.  Thus, $\pi_C(P_i)$ have
ramification index 2. The other points behave as in the previous case.  So
we have in each fiber of $\phi$ one unramified point and everything
else has ramification index 2.

\par Suppose that $P_1$ and $P_2$ are in one of the fibers which have
only one Weierstrass point.

iii) Then $w$ has to permute them, so they are not Weierstrass points.
 As in case i) $e_{\phi}\, ({\pi_C} (P_1))=4$. So there is a point of
index 4 in one of $\psi^{-1} (q_1)$, $\psi^{-1} (q_2)$,$\psi^{-1}
(q_3)$ and everything else is of ramification index 2.  The Weierstrass
points are as in case i),  unramified.

\par {\bf Case III:} Let P be the ramified point of index 3. By lemma
1, $e_\psi \, w(P)=3$.  But there is only one such point in C, so $P\in
W$. Then $e_{\pi_E \circ \psi}\,(P)=e_{\psi}\, (P) \cdot
e_{\pi_E}\, (\psi(P))=6$. So $e_{\pi_C}\, (P)\cdot e_{\phi}\, (\pi_C (P))=6$.  But
$e_{\pi_C}\, (P)=2$, because $P\in W$.  Thus, $e_{\phi}(\pi_C (P))=3$.

i) Q is in the fiber that contains three Weierstrass points. Then we
have a point of ramification index three in $\psi^{-1}(q_4)$, two
other Weierstrass points are unramified, and all the other points are
of ramification index 2.

ii) Q is in one of the fibers that  contains only one Weierstrass
point.  Then in one of $\psi^{-1} (q_1)$, $\psi^{-1} (q_2)$,  $\psi^{-1}
(q_2)$ there is a point of index 3 and everything else is of index 2.
\qed

\subsection{The case when $n$ is even}
%**************************************************************
\par Let us assume now that $deg(\psi)=n$ is an even number. The
following theorem classifies the Frey-Kani coverings in this case.

\begin{theorem}\label{thm2} 
If $n$ is an even number then the generic case 
for $\psi: C \lar E$ induce the following three cases for $\phi: \bP^1 \lar
\bP^1$:

\begin{description}
\item[I.]    $ \left ( (2)^\frac {n-2} 2  , (2)^\frac {n-2} 2  ,
(2)^\frac {n-2} 2  , (2)^ \frac {n} 2 ,  (2)  \right )
$
\item[II.]  $ \left ( (2)^\frac {n-4} 2  , (2)^\frac {n-2} 2  ,
(2)^\frac {n} 2  , (2)^ \frac {n} 2 ,  (2)             \right )
$
\item[III.] $ \left ( (2)^\frac {n-6} 2  , (2)^\frac {n} 2  ,
(2)^\frac {n} 2  , (2)^ \frac {n} 2 ,  (2)             \right )
$
\end{description}
Each of the above cases has the following degenerations (two of the branch points collapse to one)

\begin{description}
\item[I.]                                               
\begin{enumerate}
\item $\left ( (2)^\frac {n} 2  , (2)^\frac {n-2} 2  ,
 (2)^\frac {n-2} 2  , (2)^ \frac {n} 2        \right )
$
\item $\left ( (2)^\frac {n-2} 2  , (2)^\frac {n-2} 2  ,
(4) (2)^\frac {n-6} 2  , (2)^ \frac {n} 2        \right )
$
\item $\left ( (2)^\frac {n-2} 2  , (2)^\frac {n-2} 2  ,
 (2)^\frac {n-2} 2  , (4) (2)^ \frac {n-4} 2        \right )
$
\item $\left ( (3) (2)^\frac {n-4} 2  , (2)^\frac {n-2} 2  ,
 (2)^\frac {n-2} 2  ,  (2)^ \frac {n} 2        \right )
$
\end{enumerate}
\item[II.]  
\begin{enumerate}
\item	$\left ( (2)^\frac {n-2} 2  , (2)^\frac {n-2} 2  ,
(2)^\frac {n} 2  ,  (2)^ \frac {n} 2        \right )
$
\item	$\left ( (2)^\frac {n-4} 2  , (2)^\frac {n} 2  ,
(2)^\frac {n} 2  ,  (2)^ \frac {n} 2        \right )
$
\item	$\left ((4) (2)^\frac {n-8} 2  , (2)^\frac {n-2} 2  ,
(2)^\frac {n} 2  ,  (2)^ \frac {n} 2        \right )
$
\item	$\left ( (2)^\frac {n-4} 2  , (4) (2)^\frac {n-6} 2  ,
(2)^\frac {n} 2  ,   (2)^ \frac {n} 2        \right )
$
\item 	 $\left ( (2)^\frac {n-4} 2  , (2)^\frac {n-2} 2  ,
(2)^\frac {n-4} 2  ,  (2)^ \frac {n} 2        \right )
$
\item 	 $\left ((3) (2)^\frac {n-6} 2  , (2)^\frac {n-2} 2  ,
(4) (2)^\frac {n} 2  ,  (2)^ \frac {n} 2        \right )
$
\item  $\left (  (2)^\frac {n-4} 2  , (3) (2)^\frac {n-4} 2  ,
 (2)^\frac {n} 2  ,  (2)^ \frac {n} 2        \right )
$
\end{enumerate}
\item[III.] 
\begin{enumerate}
\item	 $\left ( (2)^\frac {n-4} 2  , (2)^\frac {n} 2  ,
 (2)^\frac {n} 2  , (4) (2)^ \frac {n} 2        \right )
$
\item 	 $\left ( (2)^\frac {n-6} 2  , (4) (2)^\frac {n-4} 2  ,
 (2)^\frac {n} 2  ,  (2)^ \frac {n} 2        \right )
$
\item	 $\left ( (2)^\frac {n} 2  , (2)^\frac {n} 2  ,
 (2)^\frac {n} 2  , (4) (2)^ \frac {n-10} 2        \right )
$
\item	 $\left ( (3) (2)^\frac {n-8} 2  , (2)^\frac {n} 2  ,
 (2)^\frac {n} 2  ,  (2)^ \frac {n} 2        \right )
$
\end{enumerate}
\end{description}
\end{theorem}
\proof	We know that the number of Weierstrass points in the fibers of 2-torsion points is 
$0 \mod(2)$.  Combining this with the Riemann - Hurwitz formula we get the three cases of the general case.
\par To determine the degenerate cases we consider cases when there is one branch
point for $\psi: C \lar E$.  

\par {\bf I)} First, assume that the branch point has two points $P_1$ and
$P_2$ of index $2$ (Case II).  Then $w(P_1)=P_i$ for $i=1,2$ or
$w(P_1)=P_2$. The first case implies that $P_1, P_2\in W$.  Then
$e_{\phi} (w(P_1))=e_\phi (w(P_2))=2$.  So we have case I, 1.  When
$w(P_1)=P_2$ then $e_\phi (w(P_1))=4$. Thus, we have a point of index 4
in $\phi^{-1}(q)$ for $q\in \{q_1, \dots , q_4\}$. Therefore cases 2
and 3.  If there is $P\in C$ such that $e_\psi (P)=3$, then $P\in W$
and $e_\phi (w(P))=3$. So we have case 4.

\par {\bf II)}  As in case I, if $P_1$ and $P_2$ are Weierstrass
points then they can be in the fiber of the point which has 4
 or 2 Weierstrass points. So we get two cases, namely 1 and 2.  Suppose
now that $P_1$ and $P_2$ are not Weierstrass points, thus
$w(P_1)=P_2$ and $e_\phi(w(P_1))=4$.  This point of index 4 can be in
the same fiber with 4, 2 or none Weierstrass points. So we get cases 3,
4, and 5 respectively. A point of index 3 is a Weierstrass point which
can be in the fiber which has 4 or 2 Weierstrass points. So cases 6 and 7.

\par {\bf III)}  If $P_1$ and $P_2$ are Weierstrass points then they
can be only in the fiber with 6 Weierstrass point so case 1. If they
are not then we have a point of index 4 which can be in the fiber with
all Weierstrass points or with none.  Therefore, cases 2 and 3. The
point of index 3 is a Weierstrass point so it can be in the fiber where
all the Weierstrass points are, so case 4.  This completes the proof.
\qed

\section{Maximal coverings $\psi:C \lar E$.}
%**************************************************** 

\par Let $\psi_1:C \lar E_1$ be a covering of degree $n$ from a curve of 
genus 2 to an elliptic curve.  The covering
$\psi_1:C \lar E_1$ is called a {\bf maximal covering} if it does not factor 
over a nontrivial isogeny.
A map of algebraic curves
$f: X \to Y$ induces maps between their Jacobians
$f^*: J_Y \to J_X$ and $f_*: J_X \to J_Y$.
 When $f$ is maximal then $f^*$ is injective and $ker (f_*)$ is
connected, see \cite{Se} (p. 158) for details. 

\par Let $\psi_1:C \lar E_1$ be a covering as above which is maximal.
Then ${\psi^*}_1: E_1 \to J_C$ is injective and the kernel of $\psi_{1,*}: J_C
\to E_1$ is an elliptic curve which we denote by $E_2$, see \cite{FK} or
\cite{Ku}. For a fixed Weierstrass point $P \in
C$,  we can embed $C$ to its Jacobian via
$$i_P: C \lar J_C$$
$$x \to [(x)-(P)]$$

Let $g: E_2 \to J_C$ be the natural embedding of $E_2$ in $J_C$, then
there exists $g_*: J_C \to E_2$.  Define $\psi_2=g_*\circ i_P: C \to
E_2$.
So we have the following exact sequence
$$0 \to E_2 \buildrel{g}\over\lar J_C \buildrel{\psi_{1,*}}\over\lar
E_1 \to 0$$
The dual sequence is also exact, see \cite{FK}
$$0 \to E_1 \buildrel{\psi_1^*}\over\lar J_C \buildrel{g_*}\over\lar
E_2 \to 0$$
\par The following lemma shows that $\psi_2$ has the same degree as
$\psi_1$ and is maximal. 
\begin{lemma} 
a) $deg\, (\psi_2) =n$

b) $\psi_2 $ is maximal
\end{lemma} 

\proof
\par For every $D\in Div\, (E_2)$, $deg\, (\psi_2^* D)= deg\, (\psi_2) \cdot
deg\, (D)$. Take $D={\cO }_2\in E_2$, then $deg\, (\psi_2^*\,
{\cO }_2)=deg\, (\psi_2)$. 
Also $\psi_2^*\,({\cO }_2)=(\psi_2^* \,{\cO }_2)$ as divisor and
 $$\psi_2^*\, {\cO }_2= i_P^*\, g({\cO }_2)=
i_P^*\, {\cO }_J =\psi_1^*\, {\cO }_1$$
So $deg\, (\psi_2^* \, {\cO }_2) =deg\, (\psi_1^*\, {\cO }_1)=deg\, (\psi_1)=n$

\par To prove the second part suppose $\psi_2: C \lar E_2$ is not maximal.
So there exists an elliptic curve $E_0$ and morphisms $\psi_0$ and $\beta$, such that the
 following diagram commutes

$$
\xymatrix{C \ar[d]_{\psi_0} \ar[dr]^{\psi_2}\\
E_0 \ar[r]^{\beta} & E_2\\
}
$$
Take $\psi_0 (P)$ to be the identity of $E_0$. Then exists
$\psi_{0\,*}:J_C\lar E_0$ such that $\psi_0=\psi_{0\,*}\circ i_P$. Thus,
$\psi_{2,*}=\beta\circ \psi_{0,*}$. So $ker \, \psi_{0,*}$ is a proper
subgroup of 
$ker \, \psi_{2,*}=E_1$, since $deg \, \beta > 1$.  Thus,
$$\psi_{0,*}|_{E_1}: E_1 \lar ker\, \beta$$ is a surjective
homomorphism. Therefore, $E_1$ has a proper subgroup of finite
index. So, there exists an intermediate field between function fields
$\bC (C)$ and $\bC (E_1)$. This contradicts the fact that $\psi_1$ is maximal
\qed

\par   If $deg (\psi_1)$ is an odd number then the
maximal covering $\psi_2: C \to E_2$ is unique (up to isomorphism of
elliptic curves), see Kuhn \cite{Ku}. 
\par  To each of the covers $\psi_i:C \lar E_i$, $i=1,2$, correspond Frey-Kani
covers $\phi_i: \bP^1 \lar \bP^1$.  
If the cover $\psi_1:C \lar E_1$ is given, and 
therefore $\phi_1$, we want to determine $\psi_2:C \lar E_2$ and
$\phi_2$. The study of the relation between the ramification
structures of $\phi_1$ and $\phi_2$ provides information in this direction.
  The following lemma (see \cite{FK}, p. 160)
answers this question for the set of Weierstrass points  $W=\{P_1, \dots
, P_6\}$  of C when the degree of the cover is odd. 

\par Let $\psi_i:C \lar E_i$, $i=1,2$, be maximal of odd degree
$n$. Let ${\cO}_i\in E_i[2]$ be the points which has three Weierstrass
points in its fiber.   Then we have the following: 

\begin{lemma}[Frey-Kani]  The sets $\psi_1^{-1}({\cO }_1)\cap W$ and 
$\psi_2^{-1}({\cO }_2)\cap W$ form a disjoint union of W.
\end{lemma}

\par When $n$ is even the ramification of $\psi$, is more precise.

\begin{lemma}  Let $\psi: C \lar E$ is maximal of even  degree $n$, and $Q\in E[2]$.
Then $\psi^{-1}(Q)$ has either none or two Weierstrass points. 
\end{lemma}

\proof If there are no Weierstrass points in $\psi^{-1}(Q)$ there is
nothing to prove. Suppose there is one, from lemma 3.2 we know there
are at least 2, say $P_1, P_2$. We embed $C \hookrightarrow J_C$ via
$x \lar [(x) - (P_1)]$ and $E \lar J_E$ via $x\lar [(x)-(Q)]$. 

$$
\begin{matrix}
C & \buildrel{i_{P_1}}\over\lar & J_C\\
\psi \downarrow &    & \downarrow \psi_*\\
E & \buildrel{i_Q}\over\lar & J_E
\end{matrix}
$$
Then $\psi_* ( [(x)-(P_1)])=[(\psi(x))-(Q)]$. 
\par Also, $\psi_* \psi^*=[n]$ is the multiplication by $n$ in
$E$. Since $2|n$ then $E[2]$ is a subgroup of $E[n]$. So
$\psi^*(E[2])=ker (\psi_*|_{J[2]})$, we call this group $H$. Suppose
$P_3\in \psi^{-1}(Q)$.  
Then $\psi_* (i_{P_1}(P_3))={\cO}_E$, so
$(P_1,P_3)\in H$, where the unordered pair $(P_i, P_j)$ denotes the
point $[(P_i)-(P_j)]$ of order 2 in $J_C$.
By addition of points of order 2 in $J_C$,
$(P_2,P_3)\in H$. So $H=\{0_J, (P_1,P_2), (P_1,P_3), (P_2, P_3)\}$
can't have any other points, therefore 
$\psi^{-1}(Q)$ has three Weierstrass points, which contradicts 
theorem \ref{thm2}.  Thus, there are only two Weierstrass points in
$\psi^{-1}(Q)$. 
\qed

\par The above lemma says that if $\psi$ is maximal of even degree
then the corresponding Frey-Kani covering can have only type {\bf I}
ramification,  see theorem \ref{thm1}.
%**************************************************************
\section{Arithmetic Applications}
%**************************************************************
In this section, we characterize genus 2 curves with degree 3 elliptic
subcovers and determine the j-invariants of these elliptic subcovers
in terms of coefficients of the genus 2 curve. If the elliptic
subcover is of degenerate ramification type, then its j-invariant is
determined in terms of the absolute invariants of the genus 2
curve. We find two isomorphism classes of genus 2 curves which have
both elliptic subcovers of degenerate type. 
\par When $n=5$ or 7   we discuss only Case II, iii), 
and  Case II, i)  of theorem \ref{thm1}, respectively. In both cases we
determine the j-invariants of elliptic subcovers in terms of the
coefficients of the genus 2 curves. Other types of ramifications
are computationally harder and results are very large for
display. 
%************************************************************************
\subsection{Curves of genus 2 with a degree 3 elliptic subfield.}
%***********************************************************************
Let $\psi:C\to E_1$ be a covering of degree 3, where $C$ is a  genus 2
curve given by
$$C: Y^2=x(x-1)(x-d)(x^3-ax^2+bx-c)$$
and $E_1$ an elliptic curve. 
Denote the 2-torsion points of $E_1$  by $0,1,t,s$.
Let $\phi_1$ be the Frey-Kani covering with $deg \, (\phi_1)=3$ such
that  $\phi_1\,(0)=0$, $\phi_1 \, (1)=1$, $\phi_1 \, (d)=t$, and the
roots of $f(x)=x^3-ax^2+bx-c$, are in the fiber of $s$. The fifth branch
point is infinity and in its fiber is $u$ of index 1 and infinity of
index 2. So $\phi_1$ is of generic type (Theorem \ref{thm1}).
 Points of index 2 in the fibers of 0, 1, $t$ are $m,\,n, \, p$
respectively. Then the cover is given by 
$$z=k \frac {x(x-m)^2} {x-u}$$
Then from equations: 
$$z-1=k(x-1)(x-n)^2, \quad z-t=k(x-d)(x-p)^2, \quad z-s=f(x)$$ we
 compare the coefficients of $x$ and 
 get a system of 9 equations in
the variables $a, b, c, d, k, m, n, p, t, s, u$. 
Using the Buchberger's Algorithm (see \cite{Cox}, p. 86-91)
and a computational symbolic package (as Maple) we get;
\begin{lemma} Let $E_1$ be the elliptic curve given by
$y^2=z(z-1)(z-t)(z-s)$. Then the genus 2 curve
\begin{small}
$$C: Y^2=x(x-1)\left(x-\frac {a(a-2) } {2a-3 }\right)
\left(x^3-ax^2+ \left( \frac {(2a-3)c} {(a-1)^2} +\frac {a^2} 4\right)x-c\right)$$
\end{small}
covers $E_1$ with a maximal cover of degree 3 of generic case (Theorem
\ref{thm1}).  Moreover $s$ and $t$ are given by, 
\begin{small}
$$
t = \frac {a^3(a-2)}{(2a-3)^3}, \quad s = \frac {4c}{(a-1)^2}$$
\end{small}
\end{lemma}    \qed

Next, we find the j-invariants of $E_1$ and $E_2$. The j-invariant of $E_1$ is as follows,
\begin{small}
$$j(E_1)= \frac {16} {C^2} \cdot
 \frac {A^3} {a^6 c^2 (a-1)^2 (a-2)^2(a-3)^6 \left ( (a-1)^2-4c\right )^2}$$
\end{small}
where A and C are:
\begin{small}
\begin{equation}
\begin{split}
A = & \, a^{12}-8a^{11}+16c^2a^8+11664c^2+36720c^2a^4-69984c^2a^3-192c^2a^7+77760c^2a^2\\
    & -46656c^2a + 1920c^2a^6-11232c^2a^5-4a^{10}c+26a^{10}-44a^9+41a^8-20a^7+220a^8c\\
&  -904a^7c+ 1740a^6c-1800a^5c-8a^9c-216ca^3+4a^6+972ca^4\\
C = & \, a^6-4a^5+5a^4-2a^3-32a^3c+144ca^2-216ca+108c  \\
\end{split}
\end{equation}
\end{small} 
To find $j_2$ we take  $\phi_2: \bP^1 \to \bP^1$ such that
$\phi_2(0)=\phi_2(1)=\phi_2(d)=\infty$.  Three roots of
$f_3(x)=x^3-ax^2+bx-c$ go to 2-torsion points $s_1, s_2, s_3$ of $E_2$
and 0 is the fifth branch point
of $\phi_2$.  Solving the corresponding system we get $s_1, s_2, s_3$
in terms of $a$ and $c$. Then $j_2$ is 
\begin{small}
$$j(E_2)=- \frac {16} {C} \cdot  \frac {B^3} {c\left( (a-1)^2-4c\right)}$$
\end{small}
where $A$ is as above and 
\begin{small}
$B = a^4-2a^3+a^2-24ca+36c$.
\end{small} 
%*****************************************************************
\subsection{Degenerate Cases}
%**************************************************************************
\par Notice that only one degenerate case  can occur when $n=3$.  In this case, 
one of the Weierstrass  points has ramification index 3, so the cover
is totally ramified at this point, see theorem \ref{thm1}.
\begin{lemma}
Let $E$ be an elliptic curve given by $y^2=z(z-1)(z-s)$.  Suppose
that the genus two curve $C$ with equation 
 $$Y^2=x(x-1)(x-w_1)(x-w_2)(x-w_3)$$ covers $E$, of degree 3, such
that the covering is degenerate.
Then  $w_3$ is given by

$$w_3= \frac {(4w_1^3-7w_1^2+4w_1-w_2)^3\, (4w_1^3-3w_1^2-w_2)}{16w_1^3
(w_1-1)^3\,(4w_1^3-6w_1^2+3w_1-w_2)} $$ 

and $w_1$ and $w_2$ satisfy the equation,
\begin{equation}\label{eq1}
w_1^4-4w_1^3w_2+6w_1^2w_2-4w_1w_2+w_2^2 = 0
\end{equation}
 Moreover, 
\begin{small}
$$
s= -27\left( w_1(w_1-1)\frac
{(4w_1^3-7w_1^2+4w_1-w_2)(4w_1^3-5w_1^2+2w_1-w_2)}
{(4w_1^3-9w_1^2-w_2+6w_1)(4w_1^3-3w_1^2-w_2)(4w_1^3-6w_1^2+3w_1-w_2)}
\right)^2
$$
\end{small}
\end{lemma}
\proof
We take $\psi: C\to E$ and $\phi: \bP^1\to \bP^1$ its
corresponding Frey-Kani covering.
To compute  $\phi$, let $w_1$ be the point of ramification index 3. Take
a coordinate in the lower $\bP^1$ such that $\phi\, (w_1)=0$,
$\phi\, (w_2)=s$, $\phi_2\, (w_3)=1$, and 
$\phi\, (0)=\phi\, (1)=\phi\, (\infty)=\infty$.  We denote  points of
ramification index 2 in the fibers of $s$ and 1  by $p$ and
$q$, respectively. Then, $\phi$ is given as
$z=k_2 \frac {(x-w_1)^3} {x(x-1)}$.
From the corresponding system we get the above result.
\qed

\par Denote the j-invariant of $E$ by $j_1$. 
Using the above expression of $s$ in terms of $w_1$ and $w_2$ we
get an equation in terms of $j_1$, $w_1$, and $w_2$. Taking the
resultant of this expression and equation \eqref{eq1} we get,
\begin{small}
\begin{equation}
\begin{split}
2617344w_1^2+38637j_1w_1^7-17496j_1w_1^6-29207808w_1^5-7569408w_1^3-7569408w_1^15\\
-729w_1^4j_1+5103j_1w_1^5+69984j_1w_1^9-60507j_1w_1^8+65536-589824w_1+16411392w_1^4\\
-29207808w_1^{13}+44960208w_1^{12}-60666336w_1^{11}+72010800w_1^{10}+44960208w_1^6\\
-60666336w_1^7+72010800w_1^8-75998272w_1^9+16411392w_1^{14}+2617344w_1^{16}-589824w_1^{17}\\
-60507j_1w_1^{10}+38637j_1w_1^{11}-17496j_1w_1^{12}+5103j_1w_1^{13}-729j_1w_1^{14}+65536w_1^{18}&=0
\end{split}
\end{equation}
\end{small}
\par We denote with $j$ the j-invariant of the elliptic curve
$y^2=(x-w_1)(x-w_2)(x-w_3)$. 
Then, proceeding as above, $j$ can be expressed in terms of
$w_1$ as below,
\begin{small}
\begin{equation}\label{j_1_w_1}
\begin{split}
65536w_1^6-196608w_1^5+356352w_1^4-385024w_1^3+(289536-9j)w_1^2\\
+(-129792+9j)w_1+35152-9j=0
\end{split}
\end{equation}
\end{small}
Taking the resultants of the two previous equations we have
\begin{equation}\label{eq3}
256 \, A(j)\, j_1^3 + 3\, B(j)\, j_1^2 + 6 \, C(j)\, j_1 - D(j)=0
\end{equation}
where 
\begin{small}
\begin{equation}
\begin{split}
A(j) & =(9j-35152)^4\\
B(j) & =-2187j^7+38996640j^6-277882258176j^5+998642127618048j^4\\
& -1868045010870009856j^3 +1669509508048367910912j^2\\
 & -543484034691057422696448j +16612482057244821172518912\\
C(j) &=27j^8+1125216j^7+9650655872j^6-31593875152896j^5+27748804997283840j^4\\
 & +1114515284358510673920j^3 -6061989956030939246100480j^2\\
 & +8346397859247767524611194880j +353019691006036487376293855232\\
D(j) &= (j^3+33120j^2+290490624j-310747594752)^3\\
\end{split}
\end{equation}
\end{small}
For the  genus 2 curve $C$ we compute the Igusa invariants $J_2, J_4, J_6,
J_{10}$ in terms of the coefficients of the curve,   see Igusa
\cite{Ig} for their definitions. 
The  absolute invariants of $C$ are defined it terms of Igusa
invariants as follows,
\begin{small}
\begin{equation}
i_1:=144 \frac {J_4} {J_2^2}, \quad i_2:=- 1728 \frac {J_2J_4-3J_6} {J_2^3}, \quad 
i_3 :=486 \frac {J_{10}} {J_2^5}
\label{i_u_v}
\end{equation}
\end{small}
\par Two genus 2 curves with $J_2\neq 0$ are isomorphic if and only if they have the same 
absolute invariants. 
The absolute  invariants can be expressed in terms of $w_1$ and
$w_2$. Taking the resultant of the first  two equations in \eqref{i_u_v}
we get an equation $F(i_1, i_2, w_1)=0$. The resultant of $F(i_1, i_2,
w_1)$ and equation \eqref{j_1_w_1} we get  $j =13824 \, \frac S T $ where $S$ and $T$ are:
\begin{small}
\begin{equation}
\begin{split}
S &=247945848003i_1^3-409722141024i_1^2-7591354214400i_1+17736744960000\\
 & +61379512488i_1i_2+ 64268527400i_1^2i_2-2031496516224i_2 \\
T & =
 1034723291140i_1^2i_2-3175485076512i_1i_2-7250280129792i_2+1670535171333i_1^3\\ 
& +366156782208i_1^2-67382113075200i_1+141893959680000
\end{split}
\end{equation}
\end{small}
The conjugate solutions of \eqref{eq3} are j-invariants of $E_1$ and $E_2$. 
For $j=0$ the equation (3)  has one triple root $j_1=-\frac
{1213857792}{28561}$. Then, C and E are given by,
$$Y^2=x^5-x^4+216x^2-216x$$
$$y^2=x^3-668644200   x + 6788828143125$$
For $j=1728$ the values for $j_1$ are
$$j_1=1728,\quad \frac { 942344950464}{1500625},\quad \frac { 942344950464}{1500625}$$
This value of $j$ does not give a genus 2 curve since the discriminant $J_{10}$ of $C$ is 0.

\par Next we will see what happens when both  $\phi_1$
and $\phi_2$ are degenerate. We find only two triples $(C, E_1, E_2)$
such that  the corresponding $\phi_i: C \to E_i$, $i=1,2$,  are
degenerate. It is interesting that in both cases $E_1$ and $E_2$ are
isomorphic. 
\begin{lemma}
Let  $E:y^2=z(z-1)(z-t)$ be an elliptic curve. Then the genus 2  curve 
$$Y^2=x(x-1)\left(x^3-\frac 3 2 x^2 + \frac 9 {16} x - \frac {t} {16}\right)$$
covers $E$, such that the covering is of degree 3 and the
corresponding Frey-Kani covering of  type II, iii)
(Theorem \ref{thm1}), for $t\neq 0,1$.
\end{lemma}
\proof
Let $\phi_1$ be the Frey-Kani covering with $deg \, (\phi_1)=3$ such
that $\phi_1\, (w_1)=\phi_1\, (w_2)= \phi_1\, (w_3)=t$, $\phi_1\,
(0)=0$, $\phi_1 \, (1)=1$, $\phi_1 \, (\infty)=\infty$. Let $\infty $
be the point of ramification index 3, and denote the points of
ramification index 2 in the fibers of 0 and 1 with $m$ and $n$
respectively. If $z$ is a coordinate in the lower $\bP^1$ then $\phi_1$ is given by
$z= k_1 x (x-m)^2$.
The relations  $z-1= k_1 (x-1)(x-n)^2$, $z-t=k_1 (x^3-ax^2+bx-c)$
hold, where  $x^3-ax^2+bx-c=(x-w_1)(x-w_2)(x-w_3)$. 
Comparing the  coefficients and solving the system, we get
$$(a,b,c, k_1, m, n)=\left( \frac 3 2, \frac 9 {16}, \frac t {16}, 16,
\frac 3 4 , \frac 1 4\right)$$
\qed

\par To compute  $\phi_2$, let $w_1$ be the point of ramification index 3. Take
a coordinate in the lower $\bP^1$ such that $\phi_2\, (w_1)=0$,
$\phi_2 \, (w_2)=s$, $\phi_2\, (w_3)=1$, and 
$\phi_2\, (0)=\phi_2\, (1)=\phi_2\, (\infty)=\infty$.  The points of
ramification index 2 in the fibers of $s$ and 1 we denote by $p$ and
$q$, respectively. Then $\phi_2$ is given as
$z_2=k_2 \frac {(x-w_1)^3} {x(x-1)}$.
Then from the corresponding system we get
\begin{small}
\begin{equation}
\begin{split}
& w_1 = -\frac
 {q(q-2)}{(2q-1)}\, ,
 \, w_2 = \frac {-q^3(q-2)}{(2q-1)}, \,
\, w_3 = \frac {-q(12q-8-6q^2+q^3)}{(2q-1)^3}, \\
& k_2 = \frac 1 {27} \frac {(-1+2q)^3} {q^2(q-1)^2}, \,
s = \frac {-1}{27} \frac {(-1+2q)^2(q-2)(-3q+q^3-2)}{q^2(q-1)^2}
\end{split}
\end{equation}
\end{small} 

Using the fact that the $a,b,c$ are the symmetric polynomials of
$w_1,w_2,w_3$ we have;
\begin{small}
\begin{equation}
(t,q)= \left(\frac 1 2, \frac 1 2\pm \frac 1 2 \sqrt 3\right), \, 
\left( \frac {-241+22I\sqrt 2 }{2+22I\sqrt 2}, \pm  \frac 1 2
I\sqrt2\right), \, 
\left( \frac {243} {2+22I\sqrt 2}, 1\pm \frac 1 2 I\sqrt2\right)
\end{equation}
\end{small} 
where $I=\sqrt{-1}$.
So we have three pairs of elliptic curves
\begin{small}
$$
 E_1: y^2=z(z-1)(z-\frac 1 2) \quad and \quad   E_2: y^2=z(z-1)(z+1)
$$
\end{small} 
with  $j(E_1)=j(E_2)=1728$.
\begin{small}
$$
E_1: y^2= z(z-1)\left(z- \frac    {241+22I\sqrt 2   } {-2+22I\sqrt 2}  \right), \quad
E_2: y^2= z(z-1)\left(z- \frac {241+22I\sqrt 2   } {243 }\right) 
$$
\end{small} 
with $j(E_1)=j(E_2)= \frac {-873722816} {59049}$.
\begin{small}
$$
E_1: y^2= z(z-1)\left(z- \frac {243} {1+ 2(11I\sqrt 2}\right), \quad 
E_1: y^2= z(z-1)\left(z- \frac {241 -22I\sqrt 2} {243}\right)
$$
\end{small} 
and $j(E_1)=j(E_2)= \frac {-873722816} {59049}$.
The last two cases correspond to  the same isomorphism class of
genus 2 curves.
Thus, when $\phi_1$ and $\phi_2$ are both degenerate then we get two
isomorphism classes of elliptic curves.  
Summarizing everything above we have the following table:
\begin{table}
\caption{}\label{tab1}
\renewcommand\arraystretch{1.5}
\noindent\[
\begin{array}{|c|c|c|c|c|}
\hline
f_3(x) &{E_1}&{E_2}&{j_1=j_2}\\
\hline
x^3-\frac 3 2 x^2 + \frac 9 {16} x - \frac 1 {32}      &z(z-1)(z-\frac 1 2 )&z(z-1)(z+1 )&1728\\
\hline
x^3-\frac 3 2 x^2 + \frac 9 {16}x - \frac {241+22I\sqrt{2}} {-32(1+11I\sqrt{2}) } 
&t_1= \frac {241+22I\sqrt 2} {-2+22I\sqrt 2}) 
&t_2=\frac {241+22I\sqrt 2} {243 }
&\frac {-873722816} {59049}\\
\hline
\end{array}
\]
\end{table}

where $C:Y^2=x(x-1)f_3(x)$, $E_i: y^2=z(z-1)(z-t_i)$. One can check,
using the absolute invariants of the genus two curves, that they are not
isomorphic to each other. Moreover, an equation for $E_1\iso
E_2$ in the second case is as follows:
$$y^2=z^3+z^2-277520614451197z+1880509439898307064603$$
and its conductor
$N=2^8\cdot 3\cdot 11^2\cdot 239^2\cdot 251^2$.
%************************************************************************
\subsection{Curves of genus 2 with degree 5 elliptic subfields, the 4-cycle case.}
%******************************************************************
Notice that the case II, i) does not occur when $n=5$. So we will
consider only case II, iii). We will prove the following lemma:
\begin{lemma} Let $\psi:C \to E_1$ be a covering of degree 5 such that
the corresponding  Frey-Kani cover is of ramification type II, iii)
(theorem \ref{thm1}).  Then the genus two curve can be given by
$$Y^2=x(x-1)(x-d)(x^3-ux^2+vx-w)$$
where 
\begin{small}
$$d = \frac { (3u^2-4u-4v+1)^2} {(2u-3)(6u^2-10u+5-8v)}, \quad 
w = - \frac {(u^2-6u+4v+5)(u^2-4v)}{8 (2u-3)}$$
\end{small}
 and $u$ and $v$ satisfy
\begin{small}
$$
15u^4-82u^3-8vu^2+159u^2-140u+56vu-16v^2-52v+50=0
$$
\end{small}
Moreover,  an equation of $E_1$  is $y^2=z(z-1)(z-t)$,  where 
\begin{small}
$$
t=\frac
{(u^2-4v)(-8u^4+24u^3+63u^2+64v^2-192uv+196v+16u^2v-180u+100)}{(2u-3)(6u^2-10u+5-8v) }
$$
\end{small}
\end{lemma}
\proof
Take the genus 2 curve to be
$$Y^2=x(x-1)(x-d)(x^3-ux^2+v-w)$$
Let $\phi_1$ be the Frey-Kani covering with $deg \, (\phi_1)=5$ such
that $\phi_1\, (w_1)=\phi_1\, (w_2)= \phi_1\, (w_3)=t$, $\phi_1\,
(0)=0$, $\phi_1 \, (1)=1$, and $\phi_1 \, (d)=\infty$. Take $\infty$ to be
the point of ramification index 4 such that 
$\phi_1 \,(\infty)=\infty$.  Then $\phi_1$ is given by 
$$z=k_1 \frac {x(x^2-ax+b)^2} {(x-d)}$$
Solving the corresponding system we get
the above result.  
 
\qed

From the previous lemma,  the j-invariant of the elliptic curve  satisfies 
$$F(u,v)j+G(u,v)=0$$
Taking the resultant of the previous two equations, the j-invariant
satisfies an equation of degree 2:
\begin{equation}\label{eq_j}
A(u)j^2+B(u)j+C(u)=0
\end{equation}
where 
\begin{small}
\begin{equation}
\begin{split}
A(u) & =(u-1)^2(u-2)^2(3u-4)^6(3u-5)^6(2u^2-6u+5)^8
\end{split}
\end{equation}
\begin{equation}
\begin{split}
B(u) & =-16(-7105017544704u^{33}-2816860828336128u^{31}+175917390077952u^{32}+\\
& 623116122491175945628520u^{12}+ 165647363105986609+1071822623072391493632u^{24}\\
& -697664908494919962734400u^{13}+ 10165770178171535328256u^{22}-\\
& 3521178077017962627072u^{23}- 611366039933419582356480u^{15}+ \\
& 211088208801275293447168u^{18}-117843339238828016262912u^{19}-\\
& 337258769605584067064448u^{17}+480799396622391815599360u^{16}+\\
& 58612898603387517569664u^20+139314069504u^{34}-12909484419880734720u^{27}-\\
& 284837487810868721664u^{25}+65530387559293083648u^{26}+40376325064521521748u^2-\\
& 284029170057918018876u^3-3711757861451181852u-5749828391735587589364u^5+\\
& 1452158564376272108306u^4+18345524820571264661416u^6-\\
& 48457022965012856084616u^7+ 108027612722856481764222u^8-\\
& 206208961788595840640856u^9+340743378168336968325408u^{10}-\\
& 491546319356455960291344u^{11}-25922857282984031345664u^{21}+\\
& 692593865844403162989888u^{14}+ 32784067604201472u^{30}+2146611912787372032u^{28}\\
& -295513372833693696u^{29})(2u^2-6u+5)^4\\
C(u) & =256(186624u^{16}-4478976u^{15}+50512896u^{14}-355332096u^{13}+1744993152u^{12}\\
& -6343287552u^{11}+17655393792u^{10}-38378452608u^9+65842249648u^8\\
& - 89441495616u^7+95875417216u^6-80237127456u^5 +51388251464u^4-24345314544u^3\\
& +8044840448u^2-1656421080u+160064701)^3 \\
\end{split}
\end{equation}
\end{small} 
\par The solutions of \eqref{eq_j} give the j-invariants of $E_1$ and its complement $E_2$.
\begin{example}
The two elliptic curves are isomorphic when the equation 
$$A(u)j^2+B(u)j+C(u)=0$$ of the above lemma has a double root. This happens for 
$u=\frac 3 2 \pm \frac 1 4 \sqrt{-5}$.  Then 
$$j_1=j_2=\frac {28849701763} {16941456}$$
The elliptic curve with j-invariant as above has equation, 
\begin{small}
$$y^2+yz=z^3+ 6388018241406303862z -754379181852600444980292108$$
\end{small}
\end{example}
%*************************************************
\subsection{Curves of genus 2 with degree 7 elliptic subfields, 4-cycle case}  
%***************************************************
\par The case $n=7$ is the first case that all degenerations
occur. However, it is very difficult to compute the space of genus 2
curves with degree 7 elliptic subcovers. We discuss only one
degenerate case, namely case II. iii) of theorem \ref{thm1}.
We will assume that  the genus two curve is  given by 
$$C: Y^2=x(x-1)(x-d)(x^3-ax^2+bx-c)$$ and the elliptic curve in
Legendre form $E_1: y^2=z(z-1)(z-t)$.
Moreover, let's assume that the corresponding Frey-Kani covering
$\phi:\bP^1\to \bP^1$ is of type II, i) of theorem \ref{thm1}. Take the
coordinates such that, $\phi(0)=0$, $\phi(1)=1$, $\phi(d)=t$, and
three distinct roots of $x^3-ax^2+bx-c$ are in the fiber of infinity. Let the
point of ramification index 4 be infinity, which is in the same fiber
as roots of $x^3-ax^2+bx-c$. Then the cover is given by,
$$z=k \frac {x\,P_1^2(x)} {x^3-ax^2+bx-c}$$
where $P_1(x)$ is a cubic polynomial which represents the three points
of order 2 in the fiber of 0.  Solving the corresponding system we get,
\begin{small}
\begin{equation}
\begin{split}
a  =& \frac {-1} {4A}
(7d^{20}+424t^{4}d^{8}-11072d^{12}t^{3}+2368t^{3}d^{13}-872d^{16}t^{2}-1532d^{17}t-21568d^{14}t^{2}-56d^{19}t\\  
& + 478d^{18}t+36t^{5}d-42t^{5}d^{2}+18160t^{3}d^{11}-4356t^{3}d^{10}-624t^{4}d^{6}+8t^{5}d^{3}-736t^{4}d^{7} \\
& -52594t^{2}d^{12}+624td^{14}-2576td^{15}+2725td^{16}+736td^{13}-36d^{19}-2368t^{2}d^{7}+42d^{18}\\
& + 6112d^{15}t^{2} - 29576t^{3}d^{9} - 7t^{5} +52594t^{3}d^{8} - 44496t^{3}d^{7} + 2576t^{4}d^{5}-2725t^{4}d^{4}  \\
& + 1532t^{4}d^{3} + 56t^{4}d + 872t^{3}d^{4}-6112t^{3}d^{5}-478t^{4}d^{2} - 18160d^{9}t^{2} - 424d^{12}t +
11072d^{8}t^{2} \\ 
& -8d^{17} +44496t^{2}d^{13}+ 21568t^{3}d^{6}+ 4356d^{10}t^{2}+ 29576t^{2}d^{11})\\
b = & \frac {1}{16A} 
(-14d^{21}+77d^{20}+400d^{9}t^{4}-3496t^{4}d^{8}+94280d^{12}t^{3}+1680t^{3}d^{14}-21232t^{3}d^{13}\\
&  + 1008d^{17}t^{2} + 35d^{17}t+ 31612d^{14}t^{2} + 84d^{20}t- 616d^{19}t + 1313d^{18}t - 77t^{5}d + 121t^{5}d^{2} \\ 
&-10356t^{4}d^{6}-72t^{5}d^{3}+9016t^{4}d^{7}+20t^{5}d^{4}-139344t^{2}d^{13}+269886t^{2}d^{12}-9016td^{14}\\  
 &  - 5222td^{16} + 3496td^{13} - 121d^{19} - 1680t^{2}d^{7} - 20d^{17}+ 72d^{18}+ 5352d^{15}t^{2} - 269886t^{3}d^{9}  \\
 & + 139344t^{3}d^{8}-31612t^{3}d^{7}+5222t^{4}d^{5}-35t^{4}d^{4}-5352t^{3}d^{6}-1313t^{4}d^{3}-84t^{4}d-1008t^{3}d^{4} \\ 
& + 616t^{4}d^{2} - 94280d^{9}t^{2} - 400d^{12}t + 21232d^{8}t^{2} + 219712d^{10}t^{2} - 308478t^{2}d^{11}+ 308478t^{3}d^{10} \\
& - 219712t^{3}d^{11}+ 5080t^{3}d^{5}-5080d^{16}t^{2}+10356td^{15}+ 14t^{5}        )\\
c =& -\frac {1}{448\,A}(28d^{11}-7d^{12}-561d^{4}t^{2}-1800d^{7}t+84d^{10}t+12t^{2}d+364t^{2}d^{3}-118t^{2}d^{2}+t^{3}\\
 & + 20d^{9} + 120td^{4} - 608td^{5}+ 1400td^{6} +  1311td^{8}- 42d^{10} - 140d^{6}t^{2}-504d^{9}t+ 440d^{5}t^{2})^{2}\\
\end{split}
\end{equation}
\end{small}
where,
\begin{small}
\begin{equation}
\begin{split}
A &= d(90d^{4}t^{2}-36d^{7}t-9t^{2}d-84t^{2}d^{3}+36t^{2}d^{2}+t^{3}-d^{9}+36td^{4}-90td^{5}+84td^{6}
+9td^{8}\\
& -36d^{5}t^{2})\, (168td^{6}-t^{2}-168td^{5}- 20td^{3}  + 6t^{2}d-10t^{2}d^{2} + 5t^{2}d^{3}+90td^{4}-90d^{7}t + 20td^{8} \\
&- 6d^{10} + d^{11} + 10d^{9} - 5d^{8} )  
\end{split}
\end{equation}
\end{small}
Also, $t$ and $d$  satisfy the equation,
\begin{small}
\begin{equation}
\begin{split}
 &d^{16} - 16(td^{15}+t^{3}d) + 120td^{14} - 560td^{13}+ (400t^{2} + 1420t)d^{12} - (2400t^{2} +1968t)d^{11}\\
 & + (6608t^{2} + 1400t)d^{10} - (11040t^{2}+400t)d^{9} + 12870t^{2}d^{8}- (400t^{3}+11040t^{2})d^{7} +  120t^{3}d^{2}\\
 & + (1400t^{3} + 6608t^{2})d^{6} - (1968t^{3} +2400t^{2})d^{5}+
 (1420t^{3} + 400t^{2})d^{4} - 560t^{3}d^{3}   + t^{4}=0 \\
\end{split}
\end{equation}
\end{small}
\par Thus, we can express the coefficients of $C$ in terms of $t$ and
$d$. Absolute invariants $i_1, i_2, i_3$ of $C$ can be expressed in terms
of $t$ and $d$. Using resultants and a symbolic computational package
as Maple we are able to get an equation in terms of $i_1, i_2, i_3$. 
The equation is quite large for display. 
This is the moduli space of genus two curves whose Jacobian is
the product of two elliptic curves and the Frey-Kani coverings are of
degree 7 and ramification as above. 
     
%**********************************************************                                                                  

\end{document}